\documentclass{amsart}
\usepackage{a4}
\usepackage{amssymb}

\newtheorem{Th}{Theorem}[section] \newtheorem{Cor}[Th]{Corollary}
\newtheorem{Lem}[Th]{Lemma} \newtheorem{Prop}[Th]{Proposition}

\newtheorem{Rem}[Th]{Remark}

\newtheorem{claim-num}{Claim}

\newtheorem*{lem}{Lemma}

\numberwithin{equation}{section}
\renewcommand{\theequation}{\thesection.\arabic{equation}}

\def\gl#1{\operatorname{GL}(#1)}

\def\aut#1{\operatorname{Aut}(#1)}

\def\id{\operatorname{id}}

\def\inv{^{-1}}
\def\str#1{\langle#1\rangle}

\def\f{\varphi}
\def\a{\alpha}

\def\s{\sigma}

\def\Q{\mathbf Q}

\def\av#1{\overline{#1}}

\def\cB{{\mathcal B}}

\def\cT{{\mathcal T}}
\def\cX{\mathcal X}
\def\cY{\mathcal Y}
\def\cZ{\mathcal Z}

\def\N{\mathbf N}

\def\Q{\mathbf Q}
\def\R{\mathbf R}

\renewcommand{\le}{\leqslant}
\renewcommand{\ge}{\geqslant}

\def\sym#1{\operatorname{Sym}(#1)}
\def\rank{\operatorname{rank}\,}
\def\corank{\operatorname{corank}\,}

\def\hatin{\widehat{\phantom a}}
\def\restr{\!\upharpoonright\!}
\def\invpow#1#2{{}^{#2}#1}

\begin{document}

\title
[The automorphism groups of relatively free groups]
{On Bergman's property for the automorphism groups of relatively free groups}
\author{Vladimir Tolstykh}
\date{June 28, 2004}
\address{Vladimir Tolstykh\\ Department of Mathematics\\ Yeditepe
University\\ Istanbul\\  Turkey}
\email{vtolstykh@yeditepe.edu.tr}
\subjclass[2000]{Primary: 20F28; Secondary: 20E05, 20F05}
\maketitle

\section{Introduction}

Let $\Omega$ be an infinite set and $\sym \Omega$ the
full symmetric group on $\Omega.$ In his recent
preprint \cite{Berg} Bergman proved the following delightful result.
Consider a system $(Y_i : i \in I)$ where $|I| \le |\Omega|$ of {\it subsets}
of $\sym \Omega$ whose union is $\sym \Omega.$ Then
there is a member $Y=Y_{i_0}$ of the system such that
letting $Z$ denote the conjugate set $\pi Y\pi\inv,$ where
$\pi$ is a suitable involution of $\sym \Omega$ we have that
\begin{equation}
\sym\Omega = (YZ)^4 Y \cup (ZY)^4 Z.
\end{equation}
This easily implies that the confinality of $\sym
\Omega$ is greater than $|\Omega|,$ the result first
proven by Macpherson and Neumann in \cite{MN} (the {\it confinality} of a given
group $G$ being the least cardinality
of a chain of proper subgroup whose
union is $G$.) Another consequence of the result of
Bergman's is really surprising! One of the examples
of systems $(Y_i)$ is given by the system of powers
$(X^m : m \in \N)$ of a generating set $X$ of $\sym
\Omega$ which is closed under taking inverses. By a {\it power} $X^m$ of $X$ we mean the set of
all products of the form $x_1\ldots x_m,$ where
$x_i \in X.$ Then it follows from (\theequation)
that $\sym \Omega$ equals $X^k$ for some natural
number $k.$  This naturally leads to the
following definition.

Let $G$ be a group. We say that $G$ has {\it Bergman's
property}, if for every generating set $S$ of $G$
with $S\inv=S$ we have that $G=S^k$
for some natural number $k$ \cite{DrHo,DrGo}; it might be also said
that $G$ is a {\it group of universally finite
width}. Recall that the {\it width} of a group $G$ with respect
to a generating set $S$ with $S=S\inv$ is the
least natural number $k$ such that $G=S^k,$ or
$\infty,$ otherwise.

The first example of an infinite group of universally
finite width had been found apparently by Shelah in
\cite{Sh}. The paper \cite{Sh} contains an example of
an uncountable group $G$ such that the width of $G$
relative to any generating set is at most 240.  The
striking results of Bergman's and especially his inspiring remark
in the first draft of \cite{Berg} that the further
examples of the kind might be found among `the
automorphism groups of other structures that can be
put together out of many isomorphic copies of
themselves' gave rise to the search of other examples
of groups of universally finite width. Bergman himself
suggested to analyze the situation with the
automorphism group of $\R$ as a Borel space, the
automorphism group of homogeneous Boolean spaces, the
infinite-dimensional general linear groups, the
automorphism groups of infinitely generated free
groups and some other automorphism groups.

A short while later Droste--Holland \cite{DrHo} have
proved that the automorphism group of any doubly
homogeneous chain has Bergman's property,
Droste--G\"obel \cite{DrGo} have settled Bergman's question
about the automorphism group of $\R$ as a Borel space
establishing that this group is a group of universally
finite width and the infinite-dimensional
general linear groups over divisions rings
also turned out to have Bergman's property \cite{To_GL(V)}.

Basing on (\theequation), Droste and G\"obel came to
the notion of {strong confinality} of a group. Let
$G$ be a group. Consider a {\it chain} of subsets $(Y_i : i \in
I)$ of $G$ such that

\begin{minipage}{14cm}
\begin{quote}
for all $i \in I$ the set $Y_i$ is a proper subset
of $G$ such that $Y_i\inv=Y_i$ and there exists $j=j(i) \in I$ such
that $Y_i Y_i \subseteq Y_j.$
\end{quote}
\end{minipage}
\hfill(\theequation)

Then the {\it strong confinality} of $G$ is the least
cardinality of a chain $(Y_i : i \in I)$ with
(\theequation) such that $G=\bigcup_{i\in I} Y_i.$
Bergman's results immediately imply that the strong
confinality of the full symmetric group $\sym \Omega,$
where $\Omega$ is an infinite set is greater than
$|\Omega|.$

The notions of confinality, strong confinality and
Bergman's property are embraced by the following
observation from \cite{DrHo}:  a group $G$ has
uncountable strong confinality if and only if it has
uncountable confinality and Bergman's property.
The group
$$
\text{BSym}(\Q)=\{g \in \sym \Q : \exists k \in \N \text{ s.t. }\forall x \in \Q\, |gx-x| \le k  \}
$$
of bounded permutations
of $\Q$ provides a nice illustration to this
observation: it has uncountable confinality,
but not Bergman's property \cite[Theorem 3.6]{DrGo}.

Bergman's property can be also applied for analysis of
the structural properties of groups and for answering
questions whether a given group belongs or does not
belong to a certain class of groups. For example, one
of the questions considered in
\cite{MaTho} is whether a Polish group with a comeagre
conjugacy class (the infinite symmetric
groups are examples of such groups,
by \cite{Truss}) can be written non-trivially as a free
product with amalgamation. Bergman's result provides
immediately a negative answer for $\sym \Omega$ on
an infinite set $\Omega$ (which, before
\cite{MaTho}, was known only for the
case when $\Omega$ is countable, \cite{SST}):  $\sym \Omega$
is a group of universally finite width, whereas any nontrivial
free product $G_1 *_A G_2$ with amalgamation has
infinite width relative to the generating set $G_1 \cup G_2$

In the present paper we consider the problem of
universality of finite width for the automorphism
groups of relatively free groups. We prove that the
automorphism group $\aut N$ of any infinitely
generated free nilponent group $N$ is a group of
universally finite width (and so does any homomorphic
image of $\aut N$.) We give also a partial answer to
Bergman's question about the automorphism group $\aut
F$ of infinitely generated free group $F$:
in the case when $F$ is countable, the group
$\aut F$ is indeed a group of universally finite
width.

Suppose that $\Gamma$ is a permutation group that acts
on a set $\Omega.$ We shall use standard permutation
notation, extending it, as in \cite{Berg}, to
arbitrary subsets of $\Gamma.$ Thus, if $Y$ is a
subset of $\Gamma$ and $U$ is a subset of $\Omega,$
$Y_{(U)}$ is the set of all elements of $Y$ that fix
$U$ pointwise and $Y_{\{U\}}$ is the set of all
elements of $Y$ that fix $U$ setwise. Any notation
like $Y_{*_1,*_2}$ means the set $Y_{*_1} \cap
Y_{*_2}.$

The paper is organized as follows. In Section 1
we study the automorphisms groups
of relatively free algebras. Suppose
that $\mathbf V$ is a variety algebras
and $F$ is a free algebra of $\mathbf V$
of infinite rank. Write $\Gamma$ for
$\aut F.$ Our goal is to prove that
for any Droste-G\"obel chain $(H_i : i \in I)$
with $\Gamma =\bigcup_i H_i$ of
cardinality at most $\rank F$ there is
a member $H_{i_0}$ of $(H_i)$ such that $H_{i_0}$
contains a stabilizer $\Gamma_{(U),\{W\}},$
where $F=U*W$ is the free product of
subalgebras $U$ and $W$ that are both
isomorphic to $F$ (Lemma \ref{Stab-in-a-Power}).
Lemma \ref{Stab-in-a-Power} is then
repeatedly used in the proofs of
the results on universality of
finite width for the automorphism
group of relatively free groups
we have mentioned above.

The reader may remember that a matter of generation of
a group under consideration by a pair of
stabilizers of the form $\Gamma_{(U),\{W\}}$ is one of
the central themes in many papers on confinalities (the
small index property) of certain automorphisms groups
$\Gamma.$ With this idea in mind, we introduce the
following definition. We call a variety $\mathbf
V$ of algebras a {\it BMN-variety} if for
every free algebra $F$ of infinite rank
of $\mathbf V$ we have that $\Gamma=\aut F$
is generated by any pair of stabilizers
\begin{equation}
\Gamma_{(U_1),\{U_2*W\}} \text{ and }
\Gamma_{(U_2),\{U_1*W\}},
\end{equation}
where $F=U_1 * U_2 * W$ and free factors $U_1,U_2,W$
are all isomorphic to $F.$ The results from \cite{DiNeuTho} and
\cite{Macph} can interpreted in effect that the
variety of all sets with no structure and any variety
of vector spaces over a fixed field are BMN-varieties.
With the use of Lemma \ref{Stab-in-a-Power} it is rather easy to show that
the automorphism group $\aut F$ of a free algebra $F$
of infinite rank in a BMN-variety has Bergman's
property if and only if the width of $\aut F$ relative
to any generating set of the form (\theequation) is finite.
Also, rather simple arguments show that $\aut F$ is
generated by involutions, is perfect and has confinality
greater than $\rank F.$

In Section 2 we show that variety $\mathfrak A$ of all
abelian groups and any variety $\mathfrak N_c$ of
nilpotent groups of class $\le c$ are BMN-varieties
and that the automorphism groups of relatively
free groups in these varieties have Bergman's
property. The abelian case is based upon an impressive result by
Swan which says, roughly speaking, that the
automorphism group $\aut A$ of a free abelian group
$A$ of infinite rank has finite width relative to the
set of all automorphisms of $A$ that have `unimodular
matrices' with regard to a fixed basis of $A.$ The
nilpotent case is then proved by induction
on $c.$

Methods of Section 3 are different: here we apply a
powerful fact from a paper \cite{BrEv} by
Bryant--Evans stating that a free group $F$ of
countably infinite rank possesses automorphisms that
are generic with respect to the family of all finitely
generated free factors of $F.$ Then a modification of
the famous binary tree argument from a paper
\cite{Hodges_et_al} by
Hodges--Hodkinson--Lascar--Shelah shows that $\aut F$
has Bergman's property. The scope of application
of  methods from \cite{BrEv,Hodges_et_al} is restricted by Polish
groups, and so the general case of Bergman's question
about the automorphism groups of infinitely generated free groups remains
unsettled.

The author would like to express his gratitude to
George Bergman for his comments on the author's work
\cite{To_GL(V)} and to thank Oleg Belegradek and Toma Albu
for helpful discussions.

\section{Defining a class of well-behaved varieties
of algebras}

Let $\mathfrak L$ be a language that consists only of
functional symbols. Then an {\it $\mathfrak L$-algebra} is
a structure in the language $\mathfrak L.$ An {\it $\mathfrak
L$-variety} is a class of all $\mathfrak L$-algebras
satisfying a certain set of identities over $\mathfrak L.$
{\it Free algebras} of an $\mathfrak L$-variety, their
{\it bases} and {\it ranks} are defined in the usual way of universal
algebra \cite{Cohn}.

Let $\mathbf V$ be a variety
of algebras and $F$ a free algebra of
$\mathbf V$ of infinite rank. Write
$\Gamma$ for the automorphism group of $F.$
Similarly to \cite{Macph} we call a free
factor $U$ of $F$ {\it moietous}, if
$U$ is isomorphic to $F$ and
for some decomposition
$$
F = U * W,
$$
where $*$ denotes the free product of
algebras, the free factor $W$ is also
isomorphic to $F.$ In particular,
$$
\rank F = \rank U = \rank W = \corank U = \corank W.
$$
Assume that
$$
F =\left.\prod_{i \in I}\right.^* U_i
$$
is a decomposition of $F$ into free
factors and $\s_i$ is an automorphism
of a free factor $U_i$ ($i \in I$); then
$$
\left.\prod_{i \in I}\right.^* \s_i
$$
denotes the only automorphism $\s$
of $F$ which extends all the automorphisms
$\s_i.$

The following result, the key
result of the section, generalizes
the similar result from \cite{To_GL(V)}.

\begin{Lem} \label{Stab-in-a-Power}
Let $\mathbf V$ be a variety of algebras, $F$ a free
algebra of $\mathbf V$ of infinite rank and $\Gamma$
denote the automorphism group $\aut F.$ Assume that
$\Gamma = \bigcup_{i \in I} H_i,$
where $|I| \le \rank F$ and $(H_i : i \in I)$
is a chain of proper subsets of $\Gamma$
all closed under inverses and such that
for every $i \in I$ there is $j \in I$
with $H_i H_i \subseteq H_j.$ Then an
appropriate member $H_k$ of the
chain $(H_i)$ contains a subgroup
$\Gamma_{(U),\{W\}},$ where $U,W$ are moietous free
factors of $F$ with $F=U*W.$ \end{Lem}

\begin{proof}
First, we apply a certain
`diagonal' argument to see that there exists

\stepcounter{equation}
\begin{minipage}{14cm}
\begin{quote}
a member $H_{i_0}$ of $(H_i : i \in I)$ such that for some moietous free factors $U,W$
of $F$ with $F=U*W$ the set $(H_{i_0})_{\{U\},\{W\}}$
induces the full automorphism group
$\aut U$ on $U.$
\end{quote}
\end{minipage}
\hfill(\theequation)

Let
$$
F = \left.\prod_{i \in I}\right.^* V_i
$$
be a decomposition of $F$ into
moietous free factors.
Write
$$
V_k^\times = \left.\prod_{i \ne k}\right.^* V_i
$$
for all $k \in I.$

If for some pair $(H_k,V_j)$
$$
(H_k)_{\{V_j\},\{V_j^\times\}} \text{ induces $\aut{V_j}$ on $V_j$},
$$
we are done. Suppose
otherwise. Then, in particular, for all $i$ in $I$
$$
(H_i)_{\{V_i\},\{V_i^\times\}}
\text{ does not induce $\aut{V_i}$ on $V_i$}.
$$
Hence for each $i \in I$ we can find $\s_i \in \aut{V_i}$ such
that
\begin{quote}
$\sigma_i$ does not equal to the
restriction on $V_i$ of any element from $(H_i)_{\{V_i\},\{V_i^\times\}}$.
\end{quote}
Set
$$
\sigma = \left.\prod_{i \in I}\right.^* \sigma_i.
$$
Since $\Gamma=\aut F = \bigcup_{i\in I} H_i,$ we have $\sigma \in H_m$
for some $m \in I.$ It is clear, however, that
$$
\sigma \in (H_m)_{\{V_m\},\{V_m^\times\}}.
$$
But then the restriction of $\sigma$ on $V_m$ is $\sigma_m,$
which is impossible.

So let $H_{i_0},U,W$ satisfy the condition (\theequation).
We fix some free decomposition of $U$ into
moietous free factors
$$
U = \left.\prod_{k \in \N}\right.^* U_k
$$
and for each $k \in \N$ choose a basis $\cB_k$ of $U_k.$
Let $U_0^\times$ denote the free factor
$$
\left.\prod_{k \ne 0}\right.^* U_k
$$
Consider also a family of bijections $f_{rs} : \cB_r \to \cB_s,$
where $r,s$ with $r < s$ run over $\N.$

We prove that a member of $(H_i : i \in I)$ contains
the group $\Gamma_{\{U_0\},(U_0^\times * W)}.$
Suppose $\alpha$ is any automorphism of $U_0;$
write $\a_{(k)}$ for the automorphism
$f_{0k}\a f_{0k}\inv$ of a free factor
$U_k.$ Consider the family of
all automorphisms in $H_{i_0}$ that have
the form
\begin{equation}
\gamma=\a * \left.\prod_{n \ge 1} \right.^* \a_{(3n)}\inv
* \left.\prod_{n \ge 0} \right.^* \a_{(3n+1)}
* \left.\prod_{n \ge 0} \right.^* \id_{U_{3n+2}} *\beta
\end{equation}
where $\a$ is an (arbitrary) automorphism of $U_0,$
$\beta$ is an automorphism of $W$ and
an automorphism in (\theequation) is constructed
over the free decomposition
$$
F = U_0 * \left.\prod_{n \ge 1} \right.^* U_{3n}
* \left.\prod_{n \ge 0} \right.^* U_{3n+1}
* \left.\prod_{n \ge 0} \right.^* U_{3n+2} *W.
$$
In a more simple though less formal way
any automorphism of the form (\theequation) can
be written as
$$
\gamma=\a * (\a\inv * \a\inv * \ldots ) *
(\a * \a * \ldots ) * (\id * \id * \ldots )
* \beta,
$$
Then we have
$$
\gamma\inv=
\a\inv * (\a * \a * \ldots ) *
(\a\inv * \a\inv * \ldots ) * (\id * \id * \ldots )
* \beta\inv,
$$
There exists an automorphism $\pi$ of $F$
such that

(a) $\pi$ acts as a permutation
on $\bigcup_{r \in \N} \cB_r$
and agrees on each $\cB_r$ either
with some $f_{rs}$ or with some $f_{sr}\inv$;

(b) $\pi$ fixes $W$ pointwise;

(c) and, finally, such that
$$
\pi\gamma\inv\pi\inv =
\id * (\a * \a * \ldots ) *
(\a\inv * \a\inv * \ldots ) * (\id * \id * \ldots )
* \beta\inv,
$$

This implies that
$$
\gamma\pi\gamma\inv\pi\inv =\a * (\id * \id * \ldots ) *
(\id * \id * \ldots ) * (\id * \id \ldots ) *\id,
$$
which means that $\gamma\pi \gamma\inv \pi\inv$
is an element of $\Gamma_{\{U_0\},(U_0^\times * W)}.$
Assume now that $\pi,\pi\inv$ are both
elements of some $H_{j_0}$ in the chain $(H_i : i \in I).$ Then we have
$$
\Gamma_{\{U_0\},(U_0^\times * W)} \subseteq
H_{i_0} \pi H_{i_0} \pi\inv \subseteq (H_{i_0} H_{j_0})(H_{i_0} H_{j_0})
\subseteq H_j
$$
for a suitable $j \in I,$ and the result follows.
\end{proof}

\begin{Rem}
\em In the papers \cite{BurnsPi,JCo} the method
we have applied for the construction
of $\gamma\pi \gamma\inv \pi\inv$
is referred to as `tricks
of Whitehead and Eilenberg'.
\end{Rem}

As it has been said above pairs of stabilizers of
the form $\Gamma_{(U),\{V\}}$ may generate the full
automorphism group $\Gamma.$ The known examples include,
for instance, the infinite symmetric groups \cite{DiNeuTho}
and by the general linear groups of
infinite-dimensional vector spaces \cite{Macph}.
Thus, Macpherson proves in \cite{Macph} that
given three moietous subspaces $U_1,U_2,W$ of an
infinite-dimensional vector space $V$ with
$$
V = U_1 \oplus U_2 \oplus W,
$$
we have that the group $\Gamma=\gl V$ is generated
by $\Gamma_{(U_1),\{U_2+W_2\}}$ and $\Gamma_{(U_2),\{U_1+W\}}.$
Similarly, if
$$
\Omega = U_1 \cup U_2 \cup W
$$
is a partition of an infinite set $\Omega$ into
moieties (a subset $I$ of $\Omega$ is called a {\it
moiety}, if $|I|=|\Omega \setminus I|$), then the
group $\Gamma=\sym{\Omega}$ is generated by
$\Gamma_{(U_1),\{U_2 \cup W\}}$ and
$\Gamma_{(U_2),\{U_1 \cup W\}}$ \cite{DiNeuTho}.
Moreover, the (very short) proof of the Lemma at page 580 in
\cite{DiNeuTho} shows that the following
fact is true.

\begin{Lem} \label{Berg'sLemma}
The width of the group $\Gamma =\sym \Omega$
relative to the set
$$
\Gamma_{(U_2),\{U_1 \cup W\}} \cup \Gamma_{(U_1),\{U_2 \cup W\}}
$$
is at most $3.$
\end{Lem}

Very naturally, these remarkable results
give rise to the following definition.

{\bf Definition.} Let $\mathbf V$ be
a variety of algebras. We say that
$\mathbf V$ is a {\it BMN-variety} (Bergman-Macpherson-Neumann
variety), if for every free algebra $F$ of $\mathbf V$
of infinite rank given any decomposition
$$
F = U_1 * U_2 * W
$$
of $F$ into moietous free factors, the groups
$\Gamma_{(U_1),\{U_2*W\}}$ and
$\Gamma_{(U_2),\{U_1*W\}}$ generate
the full automorphism group $\Gamma$
of $F.$

So the variety of all sets with no structure and any
variety of vector spaces over a fixed division ring
are examples of BMN-varieties. In the next section we
shall demonstrate that the variety $\mathfrak A$ and any
variety $\mathfrak N_c$ of nilpotent groups of class $\le
c$ are also examples of BMN-varieties. We remark also
that there exists a wide class of BMN-varieties of
modules.

It turns out that with the use of Lemma
\ref{Stab-in-a-Power} we could say a good
deal about the automorphism groups
of free algebras of BMN-varieties.

\begin{Th} \label{BMN:Macph'sStabs}
Let $\mathbf V$ be a BMN-variety and $F$
a free algebra of $\mathbf V$ of infinite
rank. Then the full automorphism group
$\Gamma$ of $F$ is a group of universally
finite width if and only if the width
of $\Gamma$ relative to some {\em(}any{\em)} set of the
form
$$
\Gamma_{(U_1),\{U_2*W\}} \cup
\Gamma_{(U_2),\{U_1*W\}}
$$
where $F=U_1*U_2*W$ and $U_1,U_2,W$
are moietous free factors is finite.
\end{Th}

\begin{proof}
The necessity part is obvious. Suppose that $X$ is a
generating set of $\Gamma=\aut F$ which is closed
under inverses; by adding the identity of $\Gamma$
to $X$ we can make it sure that the system
of powers $(X^k : k \in \N)$ is a chain.
The chain $(X^k : k \in \N)$ satisfies the conditions
of Lemma \ref{Stab-in-a-Power} and then there are
moietous free factors $U,W$ of $F$ with $F=U*W$ and a
power $X^m$ such that
$$
X^m \supseteq \Delta=\Gamma_{(U),\{W\}}.
$$

Now let
$$
W= W_0 * W_1
$$
be a decomposition of $W$ into moietous
free factors. Consider an automorphism
$\rho$ of $F$ that interchanges $U$
and $W_0$ and fixes $W_1$ pointwise.
Therefore by the condition the width
of $\Gamma$ relative to the
set
$$
\Delta_1=\Delta \cup \rho \Delta \rho\inv
$$
is finite: for some $k \in \N$
$$
\Gamma \subseteq \Delta_1^k.
$$
Meanwhile, if $\rho \in X^l,$ then
$$
\Gamma \subseteq \Delta_1^k \subseteq X^{k(m+2l)},
$$
and the width of $\Gamma$ relative
to $X$ is finite, as desired.
\end{proof}

Let $G$ be a group. Recall that the {\it confinality}
of $G$ is the least cardinal
$\lambda$ such that $G$ can be
expressed as the union of a chain
of $\lambda$ proper subgroups.

\begin{Th} \label{NiceProps}
Let $\mathbf V$ be a BMN-variety,
$F$ a free algebra of $\mathbf V$
of infinite rank and $\Gamma$
the automorphism group of $F.$ Then

{\em (i)} $\Gamma$ is generated by involutions;

{\em (ii)} $\Gamma$ is perfect,
that is, it coincides with the
commutator subgroup: $\Gamma=[\Gamma,\Gamma];$

{\em (iii)} the confinality of
$\Gamma=\aut F$ is greater than
$\rank F.$
\end{Th}

\begin{Rem}
\em It is interesting to compare part (ii) of the
lemma with the following corollary of Theorem C from a
paper \cite{BrRom}: if a relatively free algebra $G$
of infinite rank has the small index property, then
$\aut G$ is perfect.  \end{Rem}

\begin{proof}[Proof of Theorem {\em \ref{NiceProps}}]
(i) Let $U_1,U_2,W$ be three moietous free
factors of $F$ such that
$$
F=U_1 * W * U_2.
$$
Assume $I$ is an index set of cardinality
$\rank F$ and let
$$
(a_i :i \in I),\quad (a_i^* : i \in I)
$$
be bases of $U_1$ and $(b_i : i \in I)$
a basis of $U_2.$ Consider two (conjugate) involutions $\pi$
and $\pi_1$ which both fix $W$ pointwise and act on the bases
$(a_i)$ and $(a_i^*)$ as follows:
\begin{align*}
&\pi a_i =b_i, \qquad \forall i \in I\\
&\pi_1 a_i^* =b_i.
\end{align*}
We have therefore
that
$$
\pi_1 \pi a_i =a_i^*\quad \forall i \in I.
$$
Let now $\a$ denote the automorphism
of $U_1$ that takes the basis $(a_i)$ onto the basis $(a_i^*).$
Suppose that for all $i \in I$
$$
\a\inv a_i^* =t_i(\av a_i^*),
$$
where $t_i$ is a reduced term in the
language of $\mathbf V$ and $\av a_i^*$
denotes a finite subset of elements
of the basis $(a_i^* : i \in I).$

We then have
$$
\pi_1 \pi b_i =\pi_1 a_i =\pi_1(\a\inv a_i^*)
=\pi_1(t_i(\av a_j^*))=t_i(\av b_j)
$$
for all $i \in I.$ We see that the action of $\pi_1\pi$
on $U_2=\str{b_i : i \in I}$ is isomorphic to the
action of $\a\inv$ on $U_1=\str{a_i^* : i \in I},$
or, informally, one can write that
$$
\pi_1 \pi = \a * \id * \a\inv.
$$

Extending the principle
of the construction of $\pi_1\pi,$ one can
represent as a product of two conjugates of $\pi$ any
automorphism of $F$ of the form
\begin{equation}
\left.\prod_{n \in \N}\right.^* \a * \id *
\left.\prod_{n \in \N}\right.^* \a\inv,
\end{equation}
where the latter product corresponds
to a decomposition of $F$ into a {\it moietous}
free factors and $\a$ is the isomorphism type of an
automorphism of one of these factors. Applying tricks of Whitehead and Eilenberg as in the proof of
Lemma \ref{Stab-in-a-Power} we obtain that any automorphism
of $F$ of the form
$$
\a * \id,
$$
corresponding to a decomposition of $F$ into
moietours free factors can be written as a product of
two automorphisms of the form (\theequation), or, in
other words, this automorphism is a product of at most
four conjugates of $\pi,$ a product of at most four
involutions.

Since $\mathbf V$ is a BMN-variety, the subgroups
$\Gamma_{(U_1),\{U_2 * W\}}$ and $\Gamma_{(U_2),\{U_1*W\}}$
generate $\Gamma.$ On the other hand, each element
in these subgroups is a product of at most four
involutions.

(ii) The proof of (i) shows that $\Gamma$
is generated by all products $\rho_1 \rho_2,$ where $\rho_k$
is a conjugate of $\pi.$ Clearly, any such product
is a product of two commutators:
$$
\rho_1 \rho_2 = \s_1 \pi \s_1\inv \s_2 \pi \s_2\inv=
\s_1 \pi \s_1\inv \pi \cdot \pi \s_2 \pi \s_2\inv=
[\s_1,\pi][\pi,\s_2]
$$
(it is also easy to see that $\pi$ itself
is a commutator.)

(iii) Suppose that $(H_i : i \in I)$ is
a chain of proper subgroups of $\Gamma,$
where $|I| \le \rank F$ and the union
of $(H_i)$ is $\Gamma.$ Then by Lemma
\ref{Stab-in-a-Power} some member $H_k$ of the chain
contains a stabilizer $\Delta=\Gamma_{(U),\{W\}},$
where $U,W$ are moietous free factors whose free
product is $F.$ By the condition $\Delta$ along with
its appropriate conjugate $\rho \Delta\rho\inv$
generates $\Gamma.$ Both sets $\Delta$ and $\rho \Delta
\rho\inv$ can be found in a member, say $H_m$ of the
chain which contains $H_k$ and $\rho.$ Hence
$H_m=\Gamma,$ a contradiction.  \end{proof}

\section{Free nilpotent groups}

\begin{Prop} \label{AbCase}
{\em (i)} The variety $\mathfrak A$ of all abelian groups is
a BMN-variety;

{\em (ii)} if $A$ is an infinitely generated
free abelian group, then the automorphism
group $\aut A$ of $A$ is a group of universally
finite width.
\end{Prop}

\begin{proof}
Let $A$ be an infinitely generated free abelian group.
Swan, using a very elegant transfinite induction
argument \cite{JCo,BurnsPi}, found a set of generators of
$\aut A$ relative to which the group $\aut A$ had
finite width. More precisely, if $\cX$ is a basis of $A$ an
automorphism $\theta$ of $A$ is called
{\it $\cX$-block-unitriangular} if there is a moiety
$\cY$ of $\cX$ such that $\theta$ fixes $\cY$ elementwise and
for all $z \in \cX \setminus \cY$
$$
\theta z \equiv z\,(\mod \str{\cY}).
$$
(in particular, the `matrix' of $\theta$
with regard to the basis $\cY \cup (\cX \setminus \cY)$
is a block-unitriangular.)

\begin{Th}[Swan] \label{Swan}
The width of $\aut A$ relative
to the set of all $\cX$-block-uni\-tri\-an\-gu\-lar
automorphisms is at most $22.$
\end{Th}

The proof can be found either in \cite[Section 2]{JCo},
or in \cite[Section 2]{BurnsPi}.

Now let $U_1,U_2,W$ be moietous direct summands
of $A$ with
$$
A = U_1 \oplus U_2 \oplus W.
$$
We show that the group $\aut A$ has
finite width relative to the union
of the stabilizers
$$
X=\Gamma_{(U_1),\{U_2+W\}} \cup
\Gamma_{(U_2),\{U_1+W\}}.
$$
Suppose that $\cB_1,\cB_2$ and $\cB_W$ are bases of
$U_1,U_2$ and $W$ respectively.  Write $\cB$ for
$\cB_1 \cup \cB_2 \cup \cB_W.$ By Swan's result
and Theorem \ref{BMN:Macph'sStabs} it suffices
to show that the lengths of all $\cB$-block-unitriangular
automorphisms relative to $X$ are
bounded from above. (The {\it length}
of an element $g$ of a group $G$
with respect to a generating set $S$
of $G$ is the smallest natural $k$
such that $g$ can be expressed a product of $k$
elements of $S\cup S\inv.$)

First, note that by Lemma \ref{Berg'sLemma} any automorphism
of $A$ that acts on $\cB$ as a permutation
has the length at most $3$ with respect
to $X.$ Take a $\cB$-block-unitriangular
automorphism $\theta.$ Suppose
that $(x_i : i \in I)$ is
a basis of $U_2.$ A suitable
conjugate $\theta'=\pi\inv \theta \pi$ by
$\pi \in \aut A$ which acts on $\cB$
as a permutation preserves all elements
of $\cB_1 \cup \cB_W,$ and
for all $i \in I$ we have that
$$
\theta' x_i = x_i + a_i + b_i
$$
where $a_i \in W$ and $b_i \in U_1.$
The proof can be now completed
as in the proof of Proposition 2.2 from \cite{Macph}.
Consider $\tau_1 \in \Gamma_{(U_1),\{U_2+W\}}$
which fixes $W$ pointwise
and such that
$$
\tau_1 x_i = x_i-a_i
$$
for all $i \in I.$ Clearly, for every $i \in I$
$$
\tau_1 \theta' x_i = x_i + b_i
$$
and $\tau_1 \theta'$ acts trivially
on $\cB_1 \cup \cB_W.$ Take, further,
an automorphism $\tau_2 \in \Gamma_{(U_2),\{U_1+W\}}$
that interchanges $U_1$ and $W.$ Then
$$
\tau_2\inv \tau_1 \theta'\tau_2 x_i = x_i+a_i'  \quad \forall i \in I
$$
and $a_i' \in W.$ As above some element
$\tau_3$ of $\Gamma_{(U_1),\{U_2+W\}}$ kills
all $a_i'$:
$$
\tau_3 \tau_2\inv \tau_1 \theta'\tau_2 =\id.
$$
Thus
$$
\theta =\tau_3\inv \tau_2\inv \tau_1\inv \pi\tau_2\pi\inv
$$
and the length of $\theta$ relative to $X$
is at most $3+3+1+3=10.$
\end{proof}

\begin{Rem} \em
(i) The upper bound for the width of $\aut A$
relative to $X$ that can be extracted
from the proof of Proposition \ref{AbCase}
is of course not sharp. To improve the
said bound the  following key lemma \cite[the proof of Theorem 2.2, part (b)]{BurnsPi}
in the proof of Swan's result can be applied
rather than this result itself.

\begin{lem}
Let $\cX$ be a basis of $A.$ Then any automorphism
$\f \in \aut A$ can be written as a product
$\f_1 \f_2 \f_3$ where each $\f_i$ fixes
pointwise a moiety of $\cX$ and, moreover,
for some moiety $\cY_1$ of $\cX$ fixed by $\f_1$
pointwise the subgroup $\str{\cX \setminus \cY_1}$
is $\f_1$-invariant.
\end{lem}
(ii) It might be verified quite easily that in fact
Swan's proof works in greater generality -- in
the context of the automorphism groups
of free modules. Indeed, no changes
in the proof is required to achieve
the same conclusion as in Theorem \ref{Swan} for the automorphism
groups of infinitely generated free
modules over a ring $R$ provided
that $R$ satisfies the following
condition:

\stepcounter{equation}
\begin{minipage}{14cm}
\begin{quote}
for every free module $M$ over $R$ any direct
complement of a free direct summand of $M$ of rank one
is a free $R$-module.
\end{quote}

\end{minipage}
\hfill(\theequation)

Consequently, a stronger result than Proposition
\ref{AbCase} is true.

\begin{Th}
Let $R$ be a ring with the property {\em (\theequation)}.
Then the variety of all $R$-modules is a BMN-variety.
The automorphism group of any infinitely
generated free $R$-module is a group
of universally finite width and satisfies
the properties {\em (i-iii)} listed
in Theorem {\em \ref{NiceProps}.}
\end{Th}
\end{Rem}

We turn now to the study of the situation
with the automorphism groups of free
nilpotent groups.

\begin{Th}
{\em (i)} Any variety $\mathfrak N_c$ of all nilpotent groups
of class at most $c$ is a BMN-variety;

{\em (ii)} if $N$ is an infinitely generated
free nilpotent group, then the automorphism
group $\aut N$ of $N$ is a group of universally
finite width and satisfies
the properties {\em (i-iii)} listed
in Theorem {\em \ref{NiceProps}}.
\end{Th}

\begin{proof}
We shall prove the theorem by induction on $c.$ Proposition \ref{AbCase}
corresponds then to the case when $c=1.$

Let $N$ be an infinitely generated free nilpotent
group of nilpotency class $c$ and
let $N_c=\gamma_c(N)$ denote the $c$th term of the lower
central series of $N.$ It is well-known that $N_c$ is
a free abelian group and $N/N_c$ is a free nilpotent
group of class $c-1$ \cite{MKS}. If $\cB$ is a basis of $N$ then
the subgroup $N_c$ is generated by basic commutators
\begin{equation}
[b_1,\ldots,b_c] = [b_1,[b_2,\ldots [b_{c-1},b_c]]\ldots]
\end{equation}
where $b_1,\ldots,b_c$ are elements of $\cX$ \cite[Section 5.3]{MKS}.

Assume that
$$
N=U_1*U_2*W,
$$
where $U_1,U_2,W$ are moietous factors of $N.$ Write
$X$ for the set
$$
\Gamma_{(U_1),\{U_2*W\}} \cup
\Gamma_{(U_2),\{U_1*W\}}.
$$

According to \cite{BrMa}, the homomorphism
$\aut F \to \aut{F/\gamma_k(F)},$ where $F$ is a free
group induced by the canonical homomorphism $F \to
F/\gamma_k(F)$ is surjective for all $k \ge 2.$ Then
the homomorphism $\widehat{\phantom a} : \aut N \to
\aut{N/N_c}$ induced by the canonical homomorphism $N
\to N/N_c$ is surjective, too. It is easy to see that
the kernel $K$ of the homomorphism $\widehat{\phantom
a}$ is abelian.

The image $\widehat X$ of $X$ generates $\aut{N/N_c}$;
since by the induction hypothesis the latter group
is a group of universally finite width, then
$$
\aut{N/N_c}=\widehat X^m
$$
for some natural number $m.$ Therefore
$$
\aut N = X^m K
$$
and it remains to prove that $K$ has finite
width relative to $X.$

As $K$ is abelian, $K= K_{(U_1)} K_{(U_2*W)}.$ We show that the lengths of all elements
in $K_{(U_1)}$ relative to $X$ are uniformly bounded
from above. A similar argument can be then applied to the
elements of the subgroup $K_{(U_2*W)},$ which will
prove the result.

Suppose that $\cX$ is a basis of $U_1$
and $\cY=\cZ \cup \cT,$ where
$\cZ$ and $\cT$ are bases of $U_2$
and $W,$ respectively. Let, further,
$$
\cY = \cY_1 \cup \ldots \cup \cY_c \cup \cY_{c+1}
$$
be a partition of $\cY$ into $(c+1)$ moieties.

Take an arbitrary $\alpha \in K_{(U_1)}$
and let  $\cX=(x_i : i \in I).$ We have
$$
\alpha x_i = x_i t_i \quad \forall i \in I
$$
where $t_i$ is in $N_c.$ We write
$t_i$ as a product of the basic commutators
of the form (\theequation) over the basis $\cX \cup \cY$ and then rearrange
them to obtain a representation
$$
t_i = t_{i1} \ldots t_{i,c+1}
$$
where $t_{ij}$ is a product of basic commutators
of the form (\theequation) over the set
$$
\cX \cup (\cY \setminus \cY_j)
$$
(since each basic commutator has only
at most $c$ occurrences of elements of $\cY,$
there must be a set $\cY_k$ in which none
of these elements is contained; this
implies that each basis commutator
lies in the subgroup $\str{\cX \cup (\cY \setminus \cY_m)}$
for an appropriate $m.$)

For every $k=1,\ldots,c+1$ define the automorphism $\a_k$
of $N$ as follows:
\begin{alignat*}3
& \a_k x_i &=& x_i t_{ik} &&\quad \forall i \in I,\\
& \a_k y   &=& y          &&\quad \forall y \in \cY
\end{alignat*}
(note that if $(b_k)$ is a basis
of a free nilpotent group $G,$ then
any system $(b_k')$ of elements
of $G$ such that $b_k'$ is
congruent $b_k$ modulo $[G,G]$ for all $k$ is
also a basis of $G$ \cite[Theorem 31.25]{HNeu}). We
have therefore that $\a=\a_1 \ldots \a_{c+1}.$

Clearly, for every $k=1,\ldots,c+1$ the automorphism $\a_k$ fixes
$\cY_k$ elementwise and preserves the subgroup
$\str{\cX \cup (\cY \setminus \cY_k)}.$ This means
that a conjugate of $\a_k$ by some automorphism $\pi$
acting on $\cX \cup \cY$ as a permutation is an
element of $\Gamma_{(U_1),\{U_2*W\}}.$ Once again
by Lemma \ref{Berg'sLemma} the length of
$\pi$ with respect to $X$ does not exceed $3.$ Hence
the length of $\a_k$ relative to $X$ is at most
$3+1+3=7.$ The length of $\a$ is therefore at most
$7(c+1).$ \end{proof}

\section{Relatively free groups of countably infinite rank}

In this section we answer affirmatively Bergman's
question about the automorphism groups of free groups
of infinite rank for the case of countably
infinite rank. The very special feature of the
automorphism groups of free groups of {\it countably}
infinite rank is that they are so-called Polish
groups, that is, topological groups whose topology is
Polish. Recall that a topological space $X$ is {\it
Polish} if it is separable and there is a compatible
metric $d$ with respect to which $X$ is a complete
metric space (see \cite{Kur}; a nice introductory
account to Polish group can be also found in
\cite{Kaye}.) Developing ideas from a paper
\cite{Hodges_et_al} by Hodges--Hodkinson--Lascar--Shelah,
Bryant and Evans \cite{BrEv} introduced the
notion of automorphisms of relatively free groups that
were generic with respect to finitely generated free
factors of these groups and proved that under some
natural conditions there were `many' such generic
automorphisms.

Consider a relatively free group $G$ of countably
infinite rank. As it has been said above $\Gamma=\aut
G$ can be viewed as a Polish group: to achieve that one
defines a basis of open neighboordhoods of the identity
to be the family of all subgroups $\Gamma_{(U)},$
where $\gamma \in \Gamma$ and $U$ is a finite subset
of $G.$ A compatible metric $d$ is constructed
as follows: if $G=\{a_n : n \in \N\}$ is
an enumeration of $G$, then $d(g,h)=0,$
if $g=h$ and $1/2^n$ otherwise, when
$n$ is the minimal natural number such
that $g a_n \ne h a_n$ or $g\inv a_n \ne h\inv a_n.$

Now let us reproduce some definitions and facts from
\cite{BrEv}. Write $\cB(G)$ for the family of all finitely
generated free factors of $G.$ We say that
a tuple $(\gamma_1,\ldots,\gamma_n)$
of elements of Cartesian power $\Gamma^n$
of $\Gamma$ is {\it $\cB(G)$-generic,} if
the following two conditions are satisfied.

(1) for all $A \in \cB(G)$ the subgroups
$\Gamma_{(B)}$ for which $A \subseteq B \in \cB(G)$
and $\gamma_i B=B,$ for $i=1,\ldots,n$
form a basis of open neighbourhoods of the
identity;

(2) if $A \in \cB(G),$ $\gamma_i A =A$ for $i=1,\ldots,n,$
$A \subseteq B \in \cB(G),$ $\beta_1,\ldots,\beta_n\in
\aut B$ and $\beta\restr A = \gamma_i\restr A$
for $i=1,\ldots,n,$ then there exists $\alpha \in \Gamma_{(A)}$
such that $\gamma^\alpha_i \restr B =\beta_i$
for $i=1,\ldots,n.$ (Here and in what follows
$\gamma^\alpha=\alpha \gamma \alpha\inv.$)

Let us equip Cartesian powers $\Gamma^n$ of $\Gamma$
with product topologies. Recall that a subset $A$ of a
topological space $X$ is called {\it comeagre} if $A$
contains a countable intersection of open dense sets.
Complements of comeagre sets are called {\it meagre};
it is easy to see that a countable union of meagre
sets is again meagre. The group $\Gamma=\aut G$ is
said to have {\it ample $\cB(G)$-generic
automorphisms} if for every natural $n$ the set of
elements of $\Gamma^n$ which are $\cB(G)$-generic is
comeagre in $\Gamma^n.$

Bryant and Evans \cite{BrEv} suggested the following
property which implies ampleness of $\cB(G)$-generic
automorphisms. Let $(y_n : n \in \N)$ be a basis of
$G.$ Then $G$ has the {\it basis confinality property}
if, for every $\alpha\in \Gamma$ and every $n\in \mathbf
N$, there exist $r\in \mathbf N$ with $r\ge n$ and
$\beta \in {\rm Aut}(\langle y_1,\ldots, y_r\rangle)$
such that $\beta y_i=\alpha y_i$ for $i=1,\ldots,n$.
Lemma 1.2 from \cite{BrEv} then says that
if $\Gamma$ has basis confinality property,
then $\Gamma$ has ample $\cB(G)$-generic automorphisms.

It is easy to see that if $G$ is (absolutely) free,
then $G$ has basis confinality property and, hence,
ample $\cB(G)$-generic automorphisms \cite[Lemma
1.4]{BrEv}. Furthermore, the basis confinality
property of free groups is inherited by so-callled {\it tame}
automorphism groups of relatively free groups
\cite[Lemma 1.7]{BrEv}.  Let $\mathfrak V$ be a
variety of group and $F$ a free group of countably
infinite rank. Consider the relatively free group
$F/\mathfrak V(F)$ of $\mathfrak V,$ where $\mathfrak
V(F)$ is the verbal subgroup of $F$ corresponding to
$\mathfrak V.$ Then the automorphism group
$\aut{F/\mathfrak V(F)}$ is called {\it tame}, if the
natural homomorphism $\aut F \to \aut{F/\mathfrak
V(F)}$ induced by the canonical homomorphism $F \to
F/\mathfrak V(F)$ is surjective.  For instance, for
the varieties $\mathfrak V$ such that the relatively
free group $F/\mathfrak V(F)$ is nilpotent, the
automorphism group $\aut{F/\mathfrak V(F)}$ is tame
(\cite{BrMa}; see \cite[Theorem 1.6]{BrEv} for the other examples).

We shall use the following notation: if $\alpha$ is an
ordinal, then ${}^\alpha2$ is the set of all functions
from $\alpha$ onto $2=\{0,1\}.$ The symbol ${}^{<
\omega}2$ denotes the set $\bigcup_{n < \omega}
{}^n2.$ Notation $\sigma < \tau,$ where $\sigma,\tau$
are functions means that $\tau$ is an extension of
$\sigma.$ If $s \in {}^n2,$ then as in
\cite{Hodges_et_al} $s\hatin0$ (resp.  $s\hatin1$) is
the extension of $s$ to the ordinal $n+1$ which takes
the value $0$ (resp. $1$) at $n.$

\begin{Th}
Let $G$ be a relatively free group
of countably infinite rank which
has ample $\cB(G)$-generic automorphisms.
Then $\aut G$ is a group of universally
finite width.
\end{Th}

\begin{proof}
Suppose, towards a contradiction, that $\Gamma=\aut G$ has
infinite width relative to a generating set $X.$ By
including the identity automorphism of $G$ into $X,$ if
necessary, we may assume that the system $(X^m)$ of
powers of $X$ forms a chain. Repeating the argument
from the sufficiency part in the proof of Theorem
\ref{BMN:Macph'sStabs} one sees that some power
$X^{k_0}$ of $X$ contains a pair of stabilizers
$$
\Gamma_{(U_1),\{U_2*W\}} \text{ and }
\Gamma_{(U_2),\{U_1*W\}},
$$
where $U_1,U_2,W$ are moietous free
factors with $G=U_1*U_2*W.$ Then by Lemma
\ref{Berg'sLemma} $X^{3k_0}$ contains the subgroup of
all automorphisms of $G$ that act as permutations on a
certain basis $\cX^*$ of $G.$

Our next task is to make necessary preparations for a
$0$-$1$ game with $\cB(G)$-generic automorphisms of
$G$ (here we use the word `game' in the literal
rather than in the formal sense, though
in fact some game-theoretic argument
may be carried out, see Remark 4.2.4 in \cite{Hodges_et_al}) The reader is referred to the proof of Theorem
5.3 in \cite{Hodges_et_al} where the game has been first
described for the so-called ample homogeneous
generic automorphisms in the automorphism
groups of $\omega$-stable $\omega$-categorical
structures and for the automorphism group
of the random graph. Later Bryant and Evans
\cite{BrEv} demonstrated that the game can be adapted
for the $\cB(G)$-generic automorphisms of $G,$
thereby proving that $\aut G$ had the small
index property and uncountable
confinality.

Let us briefly discuss one of the versions of the  original
game Bryant and Evans applied in
their paper. Let $\{a_n : n < \omega\}$
be enumeration of $G$ and $(H_n)$ a countable
chain of proper subgroups of $\aut G$ whose
union is $\aut G.$ Using induction on
$s \in {}^{< \omega}2$ one constructs
of a free factor $B_s \in \cB(G)$
and elements $\gamma_s,g_{s\hatin0},g_{s\hatin1} \in \aut G$
which satisfy the following five conditions (see \cite[p. 216]{Hodges_et_al}):
\begin{itemize}
\item[G1.] $\gamma_{\varnothing}=\id$
and if $t \in {}^{<\omega}2,$ $t \le s$ and
$t \ne \varnothing,$ then $g_t B_s = B_s;$
\item[G2.] if $s \in \invpow 2n,$ then $g_{s\hatin0} \in \Gamma_{(B_s)} \cap H_n$
and $g_{s\hatin1} \in \Gamma_{(B_s)} \setminus H_n;$
\item[G3.] if $s \in \invpow 2n$ with $n >0,$ then
the tuple $\av g_s=(g_{s\restr1},\ldots,g_{s\restr n})$
is a $\cB(G)$-generic;
\item[G4.] $(g_t)^{\gamma_s}=(g_t)^{\gamma_t}$
for all $t$ such that $\varnothing < t \le s;$
\item[G5.] if $s \in \invpow 2n$ then $\gamma_{s\hatin0}\gamma_s\inv,
\gamma_{s\hatin1}\gamma_s\inv \in \Gamma_{(a_i)}
\cap \Gamma_{(\gamma_s\inv a_i)}$
for all $i \le n.$
\end{itemize}
The induction argument relies heavily on
ampleness of $\cB(G)$-generic automorphisms, since $g_{s\hatin0}$
and $g_{s\hatin1}$ are taken actually
from sets $(C \cap \Gamma_{(B_s)}) \cap H_n$
and $(C \cap \Gamma_{(B_s)}) \setminus H_n$
respectively, where $C$ is a suitable comeagre
set of $\cB(G)$-generic automorphisms.
Changing the second condition G2, the $0$-$1$
condition of the game, one can achieve
different goals: for instance, to prove (by contradiction)
that a subgroup $H$ of $\Gamma$ of small
index is open one changes the second condition
to
\begin{itemize}
\item[G$2'.$] $g_{s\hatin0} \in \Gamma_{(B_s)} \cap H$
and $g_{s\hatin1} \in \Gamma_{(B_s)} \setminus H$
\end{itemize}

Our version of the game assumes the following
variant of the condition G$2$:
\begin{itemize}
\item[G$2''$.] if $s \in {}^n2,$ then $
g_{s\hatin0} \in \Gamma_{(B_s)} \cap X^n$ and
$g_{s\hatin1} \in \Gamma_{(B_s)}\setminus X^{3n}.$
\end{itemize}

\begin{Lem} \label{MyGame}
There is a natural $m_0$ such that for every
$m \ge m_0$ and for every finite set
$E$ of $G$

{\em (i)} $\Gamma_{(E)} \cap X^m$ is not
comeagre in $\Gamma_{(E)};$

{\em (ii)} $\Gamma_{(E)} \cap X^m$ is not
meagre;

{\em (iii)} for every comeagre set $C$ of
$\Gamma$ both sets
$$
(C \cap \Gamma_{(E)}) \cap X^m \text{ and }
(C \cap \Gamma_{(E)}) \setminus X^m
$$
are not empty.
\end{Lem}

\begin{proof}
Note that $\Gamma_{(E)}$ is a Polish group.
We shall need one result from a paper
\cite{BrRom} which is an immediate
corollary of Lemma 2.1 from \cite{BrRom}.

\begin{Lem} \label{BrRomLem}
Let $\cX$ be any basis of $G$ and $E$ a
finite subset of $\cX.$ Then any automorphism
$\a \in \Gamma$ can be written as a product
$\beta \gamma,$ where $\beta \in \Gamma_{(E)}$
and $\gamma \in \Gamma_{(V_1),\{V_2\}}$, where
$V_1$ and $V_2$ are subgroups generated
by disjoint moieties of $\cX$ whose union
is $\cX.$
\end{Lem}

(i) Suppose otherwise: let $Y_m=\Gamma_{(E)} \cap X^m$
be comeagre in $\Gamma_{(E)}$ for some finite $E
\subset G.$ Take $\gamma \in Y_m.$ Since $Y_m$ is
comeagre in $\Gamma_{(E)},$ the set $\gamma Y_m$ is
also comeagre. Any two comeagre subsets always have a
non-empty interesection, and then $\gamma Y_m \cap Y_m
\ne \varnothing$ which means that
$$
\gamma \in Y_m Y_m\inv \subseteq X^{2m},
$$
whence $\Gamma_{(E)} \subseteq X^{2m}.$ Without
loss of generality we may assume that $E$ is
a subset of the above-defined basis $\cX^*.$ Any
subgroup of the form $\Gamma_{(V_1),\{V_2\}}$
from Lemma \ref{BrRomLem} over $\cX^*$
is evidently conjugate to the subgroup
$\Gamma_{(U),\{W\}}$ by a suitable
element of $\Gamma$ that acts on
$\cX^*$ as a permutation. Thus by Lemma
\ref{BrRomLem} we have
$$
\Gamma \subseteq X^{2m} X^{3k_0} X^{k_0} X^{3k_0} =X^{2m+7k_0},
$$
a contradiction.

(ii) As $\Gamma$ is Polish and $\Gamma=\bigcup X^m,$
all but finitely many of the elements of the chain
$(X^m)$ are meagre. Take as $n_0$ the minimal $m
\ge 7k_0$ such that $X^m$ is not meagre. Similarly to (i),
Lemma \ref{BrRomLem} implies that
\begin{equation}
\Gamma=\Gamma_{(E)}X^{7k_0}=\Gamma_{(E)} X^{n_0}.
\end{equation}

The index of the open subgroup $\Gamma_{(E)}$ in
$\Gamma$ is at most $\omega.$ By (\theequation)
there is a complete system $\{z_j : j \in J\}$ where
$|J| \le \omega$ of representatives of left cosets
$\gamma \Gamma_{(E)}$ which consists only of
elements of $X^{n_0}.$
Let
$$
Z_j = \{z \in X^{n_0} : z \equiv z_j (\mod \Gamma_{(E)})  \}
$$
where $j \in J.$ Since $X^{n_0} = \bigcup_{j \in J} Z_j,$
then one of the sets $Z_j$, say $Z_{j_0}$ is not meagre.
The translate $z_{j_0}\inv Z_{j_0},$ a subset
of $X^{2n_0},$ is also not meagre. We then have
$$
z_{j_0}\inv Z_{j_0} \subseteq \Gamma_{(E)} \cap X^{2n_0}
$$
and hence $\Gamma_{(E)} \cap X^m$ is not meagre
provided that $m \ge 2n_0.$

(iii) Every comeagre set meets every
nonmeagre set, and then
$$
(C \cap \Gamma_{(E)}) \cap X^m =
C \cap (\Gamma_{(E)} \cap X^m) \ne \varnothing
$$
by (ii).

The set $C^*=C \cap \Gamma_{(E)}$ is comeagre in $\Gamma_{(E)}.$
Then by (i) the relation $C^* \subseteq X^m$ is impossible.
\end{proof}

Without loss of generality we may assume
that every power of $X$ satisfies the
condition (iii) of Lemma \ref{MyGame}. Following the original game in \cite[p.
216--217]{Hodges_et_al} the reader make check now that
the properties of the $\cB(G)$-generic automorphisms
listed above and Lemma \ref{MyGame} (iii) makes it
possible to construct elements $g_{s\hatin0}, g_{s\hatin1}$
and $\gamma_s,$ where $s \in {}^n2$ satisfying the
conditions G1,G$2''$,G3,G4,G5.

We continue by analogy with the proof of Theorem 6.1 in \cite{Hodges_et_al}
making necessary adjustments.
Suppose $\sigma \in {}^\omega2$. The condition G5
implies that the sequence $(\gamma_{\sigma\restr n})$
is Cauchy; let $\gamma_\s$ denote its limit. It then
follows from G4 that
\begin{equation}
(g_{s\hatin0})^{\gamma_\sigma \gamma_\tau\inv} =
g_{s\hatin1} \quad \forall \sigma > s\hatin0,
\forall \tau  > s\hatin1
\end{equation}
where $\sigma,\tau$ are in $\invpow 2\omega.$ \cite[pp.
216--216]{Hodges_et_al}. Suppose now that $s \in {}^n2.$
By the construction $g_{s\hatin0}$ is $X^n,$ while
$g_{s\hatin1}$ is not in $X^{3n}.$ The equation
(\theequation) implies then that $\gamma_\sigma \gamma_\tau\inv$
does not belong to $X^n,$ because the length of
this element with regard to $X$ must be greater
than $n.$ In particular, if $\sigma \ne \tau,$
then $\gamma_\sigma \ne \gamma_\tau,$ since
the length of $\gamma_\sigma \gamma_\tau\inv$
relative to $X$ is greater than $1.$

Thus there are $2^\omega$ elements of the form
$\gamma_\sigma$ and hence some power
$X^n$ of $X$ contains uncountably many of them:
$$
X^n\supseteq\{\gamma_\sigma : \sigma \in \Sigma\}
$$
where $\Sigma$ is an uncountable subset
of $\invpow 2\omega.$
There exists however $m \ge 2n$ such that
for some $\s,\tau \in \Sigma$ we have
$$
\sigma \restr m=\tau\restr m \text{ and } \sigma\restr m+1\ne \tau\restr m+1.
$$
But as we saw above it follows that $\gamma_{\sigma} \gamma_{\tau}\inv,$
a product of elements of $X^m$
and also an element of $X^{2n},$ does not belong to $X^m,$ a contradiction.

\end{proof}

\begin{Cor}
Let $F$ be a free group of countably
infinite rank. Then $\aut F$ is a
group of universally finite width.
\end{Cor}

\end{document}